\documentclass[12pt]{article}
\usepackage{graphicx}
\usepackage{epsfig}
\usepackage{latexsym}
\usepackage{amsmath, amsthm, amssymb, mathrsfs}
\usepackage[fleqn,tbtags]{mathtools}

\newcommand{\beq}{\begin{equation}}
\newcommand{\eeq}{\end{equation}}

\newcommand{\vs}{\vspace{10pt}}
\newtheorem{theorem}{Theorem}
\newtheorem{remark}{Remark}
\newtheorem{lemma}{Lemma}
\newtheorem{corollary}{Corollary}

\newtheorem{definition}{Definition}

\newtheorem{example}{Example}

\def\ba{\begin{array}}
\def\ea{\end{array}}

\def\vs{\vspace{2em}}
\def\hs{\hspace{2em}}

\def\ni{\noindent}
\newcommand{\bigO}{\ensuremath{\mathcal{O}}}

\newcommand{\E}{\operatorname{E}}

\begin{document}
\pagestyle{myheadings} \markboth{Domains of Attraction}{Domains of Attraction}\markright{Domains of Attraction}

\baselineskip24pt
\title{Domains of Attraction on Countable Alphabets\footnote{{\it MSC 2010 Subject Classifications.}
Primary 6060E, 6060F; secondary 6262G.
{\it Keywords
and phrases}. Distributions on alphabets, tail index, Turing's formula, domains of attraction.} }
\author{ Zhiyi Zhang
\\
Department of Mathematics and Statistics \\ University of North Carolina at Charlotte\\
Charlotte, NC 28223 }
\date{}

\maketitle
\begin{abstract}
For each probability distribution on a countable alphabet, a sequence of positive functionals are developed as tail indices based on Turing's perspective. By and only by the asymptotic behavior of these indices, domains of attraction for all probability distributions on the alphabet are defined. The three main domains of attraction are shown to contain distributions with thick tails, thin tails and no tails respectively, resembling in parallel the three main domains of attraction, Gumbel, Fr\'{e}chet and Weibull families, for continuous random variables on the real line. In addition to the probabilistic merits associated with the domains, the tail indices are partially motivated by the fact that there exists an unbiased estimator for every index in the sequence, which is therefore statistically observable, provided that the sample is sufficiently large.
\end{abstract}

\section{Introduction and Summary.}

Consider an alphabet with countably many letters $\mathscr{X}=\{\ell_{k};k\geq 1\}$ and an associated probability distribution $P=\{p_{k};k\geq 1\}\in \mathscr{P}$ where $\mathscr{P}$ is the class of all probability distributions on $\mathscr{X}$. Let $x_{1},\cdots,x_{n}$ be an independently and identically distributed ($iid$) random sample from $\mathscr{X}$ under  $P$. Let $\{y_{k};k\geq 1\}$ and $\{\hat{p}_{k}=y_{k}/n; k\geq 1\}$ be the observed letter frequencies and relative letter frequencies in the sample.

Before proceeding further, let us first give a little thought to possible notions of an ``extreme value'' and a ``tail'' of a distribution in the current setting, as the domains of attraction are commonly discussed in association with such notions. While such notions are not required in the mathematics of this paper, it is nevertheless comforting to have them at least on an intuitive level. Unlike an $iid$ sample of a random variable on the real line where the values are numerically ordered and therefore an extreme value is naturally defined, the letters in an alphabet do not assume numerical values nor do they admit natural ordering. It is much less clear what a reasonable notion of an extreme value should be in such a case. Here if we insist to have a notion of an extreme value associated with a sample, then perhaps such a value should be based on its rarity or unusualness with respect to the observed values in the sample. The rarest values in the sample are those with frequency one and there are most commonly many more than one such observed value in a sample. If we entertain a rarer value, it has to be those with frequency zero, {\it i.e.}, the letters in the alphabet that are not represented in the sample, which, though not in the sample, are nevertheless associated with and specified by the sample. If we anticipate that another $iid$ observation from $\mathscr{X}$, say $x_{n+1}$, is to be taken, it would be reasonable then to consider the value of $x_{n+1}$ to be extreme if $x_{n+1}$ takes a letter that is not observed in the original sample of size $n$. To fix the idea, we will subsequently use the term ``an extreme value'' to mean that a new observation $x_{n+1}$ assumes a value unseen in the sample of size $n$. Similarly we can also entertain what a notation of a tail should be on an alphabet. Whenever there is no risk of ambiguity, let us loosely refer to a subset of $\mathscr{X}$ with low probability letters as a ``tail'' in the subsequent text. In this sense, a subset of $\mathscr{X}$ with very low probability letters may be referred to as a ``distant tail'', and a distribution on a finite alphabet has essentially ``no tail''. Furthermore we note that, though there is no natural ordering among the letters in $\mathscr{X}$, there is one on the index set $\{k;k\geq 1\}$. There therefore exists a natural notion of a distribution $P=\{p_{k}\}$ having a thinner tail than that of another distribution $Q=\{q_{k}\}$, in the sense of $p_{k}\leq q_{k}$ for all $k\geq k_{0}$ for some integer $k_{0}\geq 1$, when $P$ and $Q$ share a same alphabet and are enumerated by a same index set. In such a case, we will subsequently say that $P$ has a thinner tail than $Q$ {\it in the usual sense}. Finally we note that the discussion of domains of attraction for continuous random variables very much hinges on a well-defined extreme value, which is lacking on alphabets, and the differentiability of its cumulative distribution function, which is completely non-existent due to the discrete nature of alphabets. As a result of these characteristics, or the lack of them, in the current problem concerning distributions on alphabets, a fundamentally different theoretical platform is needed to move forth.

To move forth on an intuitive note, let us adopt the notation of an out-of-sample extreme value as described above. We may then entertain the probability of $x_{n+1}$ being an extreme value, {\it i.e.}, $P(\cap_{i=1}^{n}\{X_{n+1}\neq X_{i}\})$, which is, after a few algebraic steps,
\[\begin{array}{l}
  \zeta_{1,n}=\sum_{k\geq 1}p_{k}(1-p_{k})^{n}.
  \end{array}
\]
\begin{remark}
$\zeta_{1,n}$ is a member of the family of the generalized Simpson's indices $\zeta_{u,v}$ discussed by Zhang and Zhou (2010) which plays an important role in characterizing the underlying distribution $\{p_{k}\}$ (up to a permutation on the index set) and in giving alternative representations to Shannon's entropy and R\'{e}nyi's entropy, which are well-known tail indices on an alphabet,  as discussed in Zhang (2012).
\end{remark}
Clearly $\zeta_{1,n}\rightarrow 0$ as $n\rightarrow \infty$ for any probability distribution $\{p_{k}\}$ on $\mathscr{X}$.
A multiplicatively adjusted version of $\zeta_{1,n}$ is defined below and will subsequently be referred to as the tail index.
\beq \begin{array}{l}
   t_{n}=n\zeta_{1,n}=\sum_{k\geq 1}np_{k}(1-p_{k})^{n}.
 \end{array}\label{deftn}
\eeq
\begin{remark} Suppose there are two independent $iid$ samples of the same size $n$. The tail index $t_{n}$ in (\ref{deftn}) may also be interpreted as the average number of observations in one sample that are not found in the other sample.
\end{remark}

The fact that $t_{n}$ is tail-relevant is manifested in the fact that $\zeta_{1,n}$ is tail-relevant. To see that $\zeta_{1,n}$ is tail-relevant, let us first consider $\pi_{0}=\sum_{k\geq 1}p_{k}1[y_{k}=0]$. $1-\pi_{0}$ is often referred to as the sample coverage of a population in the literature. Since the letters not represented in a large sample are likely those with low probabilities, it is reasonable to think that $\pi_{0}$ is a tail-relevant quantity for a large $n$; and yet
$\zeta_{1,n}=\E(\pi_{0})$.
Intuitively one would expect $\pi_{0}$ to take a smaller (larger) value under a more (less) concentrated probability distribution, and therefore to expect $\zeta_{1,n}$, and hence $t_{n}$, to be a reasonable measure to characterize the tail of a distribution on an alphabet. Also to be noted is that, for any given integer $k_{0}\geq 1$, the first $k_{0}$ terms in the re-expression of $t_{n}$ below converges to zero exponentially fast as $n\rightarrow \infty$
\[ \begin{array}{l}t_{n}=\sum_{k\leq k_{0}}np_{k}(1-p_{k})^{n}+\sum_{k>k_{0}}np_{k}(1-p_{k})^{n},\end{array}
\] and therefore the asymptotic behavior of $t_{n}$ has essentially nothing to do with how the probabilities are distributed over any fixed and finite subset of $\mathscr{X}$, further noting that $t_{n}$ is invariant under any permutation on the index set $\{k\}$.

\begin{remark}
 Good (1953) introduced a remarkable estimator of $\pi_{0}$ in the form of $N_{1}/n$ where $N_{1}=\sum_{k}1[y_{k}=1]$. The estimator, also known as Turing's formula, is the subject of much research in the existing literature. Notable papers on this topic include Robbins (1968) and Esty (1983), and more recent advances are reported in Zhang and Huang (2008), Zhang and Zhang (2009) and Zhang (2013). One of the most intriguing characterisitcs of Turing's formula is its ability to infer nonparametrically the probability beyond the range of observed data.
\end{remark}
\begin{remark}
 Domains of attraction for distributions of continuous random variables are a long-standing focal point of the extreme value theory. The large volume of research on this topic in the existing literature goes back to Fr\'{e}chet (1927) and Fisher and Tippett (1928), and includes full analyses by Gnedenko (1944) and Smirnov (1949). There the three main domains of attraction are defined along the lines of Gumbel family (thick tails), Fr\'{e}chet family (thin tails) and Weibull family (no tails). The main objective of this paper is to similarly characterize many distributions on alphabets by the indices $\{t_{n},n\geq 1\}$ into three domains, Domain 0 (no tails), Domain 1 (thin tails), and Domain 2 (thick tails).
\end{remark}

\begin{definition}
A distribution $P=\{p_{k}\}$ on $\mathscr{X}$ is said to belong to
  \begin{enumerate}
    \item Domain 0 if $\lim_{n\rightarrow \infty}t_{n}=0$,
    \item Domain 1 if $\limsup_{n\rightarrow \infty} t_{n}=c_{P}$ for some constant $c_{P}>0$,
    \item Domain 2 if $\lim_{n\rightarrow \infty}t_{n}=\infty$, and
    \item Domain $T$, or Domain Transient, if it does not belong to Domains 0, 1, or 2.
  \end{enumerate}
\label{def1}
\end{definition}

The four domains so defined above form a partition of $\mathscr{P}$. The primary results established in this paper include:
\begin{enumerate}
  \item Domain 0 does and only does include probability distributions with positive probabilities on a finite subset of $\mathscr{X}$.
  \item Domain 1 includes distributions with thin tails such as $p_{k}=\bigO \left(a^{-\lambda k}\right)$, $p_{k}=\bigO \left(a^{-\lambda k^{2}}\right)$, and $p_{k}=\bigO \left(k^{r}a^{-\lambda k}\right)$ where $a>1$, $\lambda>0$ and $r\in (-\infty, \infty)$ are constants.
  \item Domain 2 includes distributions with thick tails such as $p_{k}=\bigO \left(k^{-\lambda}\right)$ and $p_{k}=\bigO \left((k\ln^{\lambda} k)^{-1}\right)$ where $\lambda>1$.
  \item A relative regularity condition between two distributions (one dominates the other) is defined. Under this condition, all distributions on a countably infinite alphabet, that are dominated by a Domain 1 distribution, must also belong to Domain 1.
  \item Domain $T$ is not empty.
\end{enumerate}
The secondary results established in this paper include:
\begin{enumerate}
  \item In Domain 0, $t_{n}\rightarrow 0$ exponentially fast for every distribution.
  \item The tail index $t_{n}$ of a distribution with tail $p_{k}=\bigO \left(e^{-\lambda k}\right)$ where $\lambda>0$ in Domain 1 perpetually oscillates between two positive constants and does not have a limit as $n\rightarrow \infty$.
  \item There is a uniform positive lower bound for $\limsup_{n\rightarrow \infty}t_{n}$ for all distributions with positive probabilities on infinitely many letters of $\mathscr{X}$.
\end{enumerate}

All above mentioned results are given in Section 2. Section 3 includes several constructed examples, each of which illustrate a point of interest. The paper ends with a brief discussion in Section 4 on the statistical implication of the established results.

\section{Main Results.}
Let $K$ be the effective cardinality, or simply the cardinality when there is no ambiguity, of $\mathscr{X}$, {\it i.e.}, $K=\sum_{k}1[p_{k}>0]$.
\begin{lemma}
  If $K=\infty$, then there exists a subsequence $\{n_{k};k\geq 1\}$ in $\mathbb{N}$, satisfying $n_{k}\rightarrow \infty$ as $k\rightarrow \infty$, such that $t_{n_{k}}>c>0$ for all sufficiently large $k$.
\label{lemma1}
\end{lemma}\ni {\it Proof}.
Let us assume without loss of generality that $p_{k}>0$ for all $k\geq 1$.
Since $\zeta_{1,n}$ is invariant with respect to any permutation on the index set $\{k;k\geq 1\}$, it can be assumed without loss of generality that $\{p_{k}\}$ is non-increasing in $k$. For every $k$, let $n_{k}=\lfloor 1/p_{k}\rfloor$. With $n_{k}$ so defined, we have $1/(n_{k}+1)<p_{k}\leq 1/n_{k}$ for every $k$ and $\lim_{k\rightarrow \infty}n_{k}=\infty$ though $\{n_{k}\}$ may not necessarily be strictly increasing. By construction, the following are true about the $n_{k}$, $k\geq 1$.
 \begin{enumerate}
   \item $\{n_{k};k\geq 1\}$ is an infinite subset of $\mathbb{N}$.
   \item Every $p_{k}$ is covered by the interval $(1/(n_{k}+1),1/n_{k}]$.
   \item Every interval $(1/(n_{k}+1),1/n_{k}]$ covers at least one $p_{k}$ and at most finitely many $p_{k}$s.
 \end{enumerate}

Let $f_{n}(x)=nx(1-x)^{n}$ for $x\in [0,1]$. $f_{n}(x)$ attains its maximum at $x=(n+1)^{-1}$ with value
\[\begin{array}{l}
    f_{n}\left(\frac{1}{n+1}\right)=\frac{n}{n+1}\left(1-\frac{1}{n+1}\right)^{n}=\left(\frac{n}{n+1}\right)^{n+1}\rightarrow e^{-1}.
  \end{array}\] Also we have
\[\begin{array}{l}
    f_{n}\left(\frac{1}{n}\right)=\left(1-\frac{1}{n}\right)^{n}\rightarrow e^{-1}.
  \end{array}
\] Furthermore since $f'(x)<0$ for $(n+1)^{-1}<x<1$, we have
\[\begin{array}{l}
    f_{n}\left(\frac{1}{n}\right)<f_{n}(x)<f_{n}\left(\frac{1}{n+1}\right)\hs \mbox{for} \hs \frac{1}{n+1}<x<\frac{1}{n}.
  \end{array}
\]

Since $f_{n}(1/n)\rightarrow e^{-1}$ and $f_{n}(1/(n+1))\rightarrow e^{-1}$, for any arbitrarily small but fixed $\varepsilon>0$ there exists a positive $N_{\varepsilon}$ such that for any $n>N_{\varepsilon}$, $f_{n}(1/(n+1))>f_{n}(1/n)>e^{-1}-\varepsilon$.

Since $\lim_{k\rightarrow \infty}n_{k}=\infty$ and $\{n_{k}\}$ is non-decreasing, there exists an integer $K_{\varepsilon}>0$ such that $n_{k}>N_{\varepsilon}$ for all $k>K_{\varepsilon}$. Consider the sub-sequence $\{t_{n_{k}}; k\geq 1\}$. For any $k>K_{\varepsilon}$,
\[\begin{array}{l}
    t_{n_{k}}=\sum_{i=1}^{\infty}n_{k}p_{i}(1-p_{i})^{n_{k}}>f_{n_{k}}(p_{k}).
  \end{array}
\] Since $p_{k}\in (1/(n_{k}+1),1/n_{k}]$ and $f_{n_{k}}(x)$ is decreasing on the interval $(1/(n_{k}+1),1/n_{k}]$, we have
\[\begin{array}{l}
    f_{n_{k}}(p_{k})>f_{n_{k}}\left(\frac{1}{n_{k}}\right)\geq e^{-1}-\varepsilon,
  \end{array}
\] and hence $t_{n_{k}}>f_{n_{k}}(p_{k})\geq e^{-1}-\varepsilon$ for all $k>K_{\varepsilon}$. \hfill $\Box$

\begin{theorem} $K<\infty$ if and only if
\beq \lim_{n\rightarrow \infty}t_{n}=0.\label{condition1}\eeq
\label{theorem1}
\end{theorem}\ni {\it Proof}. Assuming that $P=\{p_{k}; 1\leq k\leq K\}$ where $K$ is finite and $p_{k}>0$ for all $k$, $1\leq k\leq K$, and denoting $p_{0}=\min\{p_{k};1\leq k\leq K\}>0$, the necessity of (\ref{condition1}) follows the fact that as $n\rightarrow \infty$
\[ \begin{array}{l}
    t_{n}=n\sum_{k}^{K}p_{k}(1-p_{k})^{n}\leq n\sum_{k}^{K}p_{k}(1-p_{0})^{n}=n(1-p_{0})^{n}\rightarrow 0.
   \end{array}
\]The sufficiency of (\ref{condition1}) follows the fact that, if $K=\infty$, then Lemma \ref{lemma1} would provide a contradiction to (\ref{condition1}).
\hfill $\Box$

In fact the proof of Theorem \ref{theorem1} also establishes the following corollary.
\begin{corollary}
 $K<\infty$ if and only if $t_{n}\leq \bigO(nq_{0}^n)$ where $q_{0}$ is a constant in $(0,1)$.
\label{coro1}
\end{corollary}

Theorem \ref{theorem1} and Corollary \ref{coro1} firmly characterize Domain 0 as a family of distributions on finite alphabets. All distributions outside of Domain 0 must have positive probabilities on infinitely many letters of $\mathscr{X}$. The entire class of such distributions is denoted as $\mathscr{P}_{+}$. In fact in the subsequent text when there is no ambiguity $\mathscr{P}_{+}$ will denote the entire class of distributions with a positive probability on every $\ell_{k}$ in $\mathscr{X}$. For all distributions in $\mathscr{P}_{+}$, a natural group would be those for which $\lim_{n}t_{n}=\infty$ and so Domain 2 is defined.

The following three lemmas are useful in the proof of Theorem \ref{theorem2} below which puts distributions with a power decaying or a slower tail in Domain 2. Lemma
\ref{euler} is a version of the well-known Euler-Maclaurin formula
and therefore is referred to as the Euler-Maclaurin Lemma subsequently.
\begin{lemma} (Euler-Maclaurin) Let $f_{n}(x)$ be a continuous function of $x$ on
$[x_{0}, \infty)$ where $x_{0}$ is a positive integer. Suppose
$f_{n}(x)$ is increasing on $[x_{0},x(n)]$ and decreasing on
$[x(n),\infty)$. If $f_{n}(x_{0})\rightarrow 0$ and
$f_{n}(x(n))\rightarrow 0$, then
\[\begin{array}{l}\lim_{n\rightarrow \infty} \sum_{k\geq
x_{0}}f_{n}(k)=\lim_{n\rightarrow \infty}\int_{x_{0}}^{\infty}f_{n}(x)dx. \end{array}\]
\label{euler}
\end{lemma}
\noindent {\it Proof.} It can be verified that
\[\begin{array}{l}
 \sum_{x_{0}\leq k\leq x(n)}f_{n}(k)-f_{n}(x(n))\leq
 \int_{x_{0}}^{x(n)}f_{n}(x)dx\leq
 \sum_{x_{0}+1\leq k<x(n)}f_{n}(k)+f_{n}(x(n))\hspace{1em}and\\ \\
  \sum_{k>x(n)}f_{n}(k)-f_{n}(x(n))\leq
 \int_{x(n)}^{\infty}f_{n}(x)dx\leq
 \sum_{k\geq x(n)}f_{n}(k)+f_{n}(x(n)).
\end{array}
\]Adding the corresponding parts of the two expressions above and
taking limits give
\[\begin{array}{l}
 \lim_{n\rightarrow \infty} \sum_{k=x_{0}}^{\infty}f_{n}(k)-2\lim_{n\rightarrow \infty}
 f_{n}(x(n))\leq \lim_{n\rightarrow \infty}
 \int_{x_{0}}^{\infty}f_{n}(x)dx \\ \\
 \leq \lim_{n\rightarrow \infty} \sum_{k=x_{0}}^{\infty}f_{n}(k)-\lim_{n\rightarrow \infty} f_{n}(x_{0})+2\lim_{n\rightarrow \infty} f_{n}(x(n)).
 \end{array}
\]The desired result follows the conditions of the lemma.\hfill
$\Box$

The next lemma includes two trivial but useful facts.
\begin{lemma}
  \begin{enumerate}
    \item For any real number $x\in [0,1)$, $ 1-x\geq \exp\left(-\frac{x}{1-x}\right).$
    \item For any real number $x\in (0,1/2)$, $ \frac{1}{1-x}<1+2x.$
  \end{enumerate}
  \label{2facts}
\end{lemma}
\ni {\it Proof}. 
For part 1, the function
$y=\frac{1}{1+t}e^{t}$ is strictly increasing over $[0,\infty)$,
and has value $1$ at $t=0$. Therefore $\frac{1}{1+t}e^{t}\geq 1 $
for $t\in [0,\infty)$. The desired inequality follows the change
of variable $x=t/(1+t)$. For part 2, the proof is trivial. \hfill $\Box$

\begin{lemma} For any given probability distribution $P=\{p_{k}; k\geq 1\}$,  $n^{1-\delta}\sum_{k}p_{k}(1-p_{k})^{n}\rightarrow c>0$ for some constants $c>0$ and $\delta \in (0,1)$,
if and only if $n^{1-\delta}\sum_{k}p_{k}e^{-np_{k}}\rightarrow c>0$,  as $n\rightarrow \infty$.
\label{lemma2}
\end{lemma}
\ni {\it Proof}. Let $\delta^{*}=\delta/8$. Consider the partition of the index set $\{k;k\geq 1\}=I\cup II$ where
\[\begin{array}{l}
   I=\{k;p_{k}\leq 1/n^{1-\delta^{*}}\}\hs \mbox{and}\hs II=\{k;p_{k}> 1/n^{1-\delta^{*}}\}.
   \end{array}
\] Since $pe^{-np}$ has a
negative derivative with respect to $p$ on interval $(1/n, 1]$ and
hence on $(1/n^{1-\delta^{*}},1]$ for large $n$, $p_{k}e^{-np_{k}}$
attains its maximum at $p_{k}=1/n^{1-\delta^{*}}$ for every $k\in II$. Therefore noting
that there are at most $n^{1-\delta^{*}}$ indices in $II$,

\[\begin{array}{l}
   0\leq n^{1-\delta}\sum_{II}p_{k}(1-p_{k})^{n}\leq n^{1-\delta}\sum_{II}p_{k}e^{-np_{k}}  \\ \\
   \hspace{2em}\leq n^{1-\delta}\sum_{II}\left(\frac{1}{n^{1-\delta^{*}}}e^{-\frac{n}{n^{1-\delta^{*}}}}
   \right)  \leq n^{1-\delta}n^{1-\delta^{*}}\left(\frac{1}{n^{1-\delta^{*}}}e^{-\frac{n}{n^{1-\delta^{*}}}}\right)\\ \\
   \hspace{2em}
   =n^{1-\delta}e^{-n^{\delta^{*}}}
   \rightarrow 0.
  \end{array}
\]

Thus
 \beq \begin{array}{c}\lim_{n\rightarrow\infty}n^{1-\delta}\sum_{k}p_{k}(1-p_{k})^{n}= \lim_{n\rightarrow\infty} n^{1-\delta}\sum_{I}p_{k}(1-p_{k})^{n}
 \end{array}
 \label{cc1}
 \eeq and
 \beq \begin{array}{l}
  \lim_{n\rightarrow\infty}  n^{1-\delta}\sum_{k} p_{k}e^{-np_{k}}= \lim_{n\rightarrow\infty} n^{1-\delta}\sum_{I}p_{k}e^{-np_{k}}.
  \end{array}
  \label{cc2}
  \eeq

On the other hand, since $1-p\leq e^{-p}$ for all $p\in [0,1]$,
 \[ \begin{array}{l}
   n^{1-\delta}\sum_{I}p_{k}(1-p_{k})^{n}\leq n^{1-\delta}\sum_{I}p_{k}e^{-np_{k}}.
  \end{array}
 \]

Furthermore, applying 1) and 2) of Lemma \ref{2facts} in
the first and the third steps below respectively leads to
 \[ \begin{array}{l}
   n^{1-\delta}\sum_{I}p_{k}(1-p_{k})^{n}\geq
   n^{1-\delta}\sum_{I}p_{k}\exp\left(-\frac{np_{k}}{1-p_{k}}\right)\\
   \\ \hspace{2em}
   \geq n^{1-\delta}\sum_{I}p_{k}\exp\left(-\frac{np_{k}}{1-\sup_{I}p_{k}}\right)  \geq n^{1-\delta}\sum_{I}\exp(-2n(\sup_{I}p_{k})^{2})p_{k}e^{-np_{k}}.
  \end{array}
  \]

 Noting the fact that $\lim_{n\rightarrow\infty}
 \exp(-2n(\sup_{I}p_{k})^{2})=1$ uniformly by the definition of $I$,
 \[\begin{array}{l}
   \lim_{n\rightarrow\infty} n^{1-\delta}\sum_{I}p_{k}(1-p_{k})^{n}= \lim_{n\rightarrow\infty}
 n^{1-\delta}\sum_{I}p_{k}e^{-np_{k}},
 \end{array}\]and hence, by
 (\ref{cc1}) and (\ref{cc2}), the lemma follows. \hfill $\Box$

\begin{theorem}
 For any given probability distribution $P=\{p_{k}; k\geq 1\}$, if there exists constants $\lambda>1$, $c>0$ and integer $k_{0}\geq 1$ such that for all $k\geq k_{0}$
  \beq p_{k}\geq ck^{-\lambda},\label{condition3}\eeq then $\lim_{n\rightarrow \infty}t_{n}=\infty$.
\label{theorem2}
\end{theorem}\ni {\it Proof}. For clarity, the proof is given in 2 cases respectively:
\begin{enumerate}
  \item $p_{k}= ck^{-\lambda}$ for all $k\geq k_{0}$ for some $k_{0}>1$, and
  \item $p_{k}\geq ck^{-\lambda}$ for all $k\geq k_{0}$ for some $k_{0}>1$.
\end{enumerate}
\ni {\it Case 1}: Assuming $p_{k}= ck^{-\lambda}$ for all $k\geq k_{0}$, it suffices to consider the partial series
$\sum_{k\geq k_{0}}np_{k}(1-p_{k})^{n}$. First consider
\[ \begin{array}{l}
     n^{1-\frac{1}{\lambda}}\sum_{k=k_{0}}^{\infty}p_{k}e^{-np_{k}}=n^{1-\frac{1}{\lambda}}\sum_{k=k_{0}}^{\infty}ck^{-\lambda}e^{-nck^{-\lambda}}
     =\sum_{k=k_{0}}^{\infty}f_{n}(k)
   \end{array}
\] where $f_{n}(x)=n^{1-\frac{1}{\lambda}}cx^{-\lambda}e^{-ncx^{-\lambda}}$. Since it is easily verified that
\[ \begin{array}{l}
   f'_{n}(x)=-\lambda c n^{1-\frac{1}{\lambda}}x^{-(\lambda+1)}(1-ncx^{-\lambda})e^{-ncx^{-\lambda}},
\end{array}
\] it can be seen that, $f_{n}(x)$ increases over $[1, (nc)^{1/\lambda}]$ and decreases over $[(nc)^{1/\lambda}, \infty)$. Let
$x_{0}=k_{0}$   and $x(n)=(nc)^{1/\lambda}$. It is clear that $f_{n}(x_{0})\rightarrow 0$ and
\[
f_{n}(x(n))=n^{1-\frac{1}{\lambda}}c(nc)^{-1}e^{-nc(nc)^{-1}}=n^{1-\frac{1}{\lambda}}c(nc)^{-1}e^{-1}=\frac{1}{en^{1/\lambda}}\rightarrow 0.
\] Invoking the Euler-Maclaurin Lemma, we have, with changes of variable $t=x^{-\lambda}$ and then $s=nct$,
\[\begin{array}{l}
   n^{1-\frac{1}{\lambda}}\sum_{k=k_{0}}^{\infty}p_{k}e^{-np_{k}}\sim \int_{x_{0}}^{\infty}n^{1-\frac{1}{\lambda}}cx^{-\lambda}e^{-ncx^{-\lambda}}dx
   =\frac{c}{\lambda}\int_{0}^{x_{0}^{-\lambda}}n^{1-\frac{1}{\lambda}}t^{-\frac{1}{\lambda}}e^{-nct}dt \\ \\
 \hs =\frac{c}{\lambda}n^{1-\frac{1}{\lambda}}\int_{0}^{x_{0}^{-\lambda}}(nct)^{-\frac{1}{\lambda}}(nc)^{-1+\frac{1}{\lambda}}e^{-nct}d(nct)
     =\frac{c}{\lambda}n^{1-\frac{1}{\lambda}}(nc)^{-1+\frac{1}{\lambda}}\int_{0}^{ncx_{0}^{-\lambda}}s^{-\frac{1}{\lambda}}e^{-s}ds \\ \\
 \hs =\frac{c^{\frac{1}{\lambda}}}{\lambda }n^{0}\int_{0}^{ncx_{0}^{-\lambda}}s^{-\frac{1}{\lambda}}e^{-s}ds
     =\frac{c^{\frac{1}{\lambda}}}{\lambda }\int_{0}^{ncx_{0}^{-\lambda}}s^{\left(1-\frac{1}{\lambda}\right)-1}e^{-s}ds \\ \\
 \hs =\frac{c^{\frac{1}{\lambda}}}{\lambda }\Gamma\left(1-\frac{1}{\lambda}\right)
     \left[\frac{1}{\Gamma\left(1-\frac{1}{\lambda}\right)}\int_{0}^{ncx_{0}^{-\lambda}}s^{\left(1-\frac{1}{\lambda}\right)-1}e^{-s}ds\right]
     \rightarrow \frac{c^{\frac{1}{\lambda}}}{\lambda }\Gamma\left(1-\frac{1}{\lambda}\right)>0.
\end{array}
\] Hence by Lemma \ref{lemma2},
$n^{1-1/\lambda}\sum_{k=1}^{\infty}p_{k}(1-p_{k})^{n}\rightarrow
c^{1/\lambda}\lambda^{-1}\Gamma\left(1-1/\lambda\right)>0$ and therefore $t_{n}\rightarrow \infty$.

\vs
\ni {\it Case 2}: Assuming $p_{k}\geq ck^{-\lambda}=\vcentcolon q_{k}$ for all $k\geq k_{0}$ for some $k_{0}\geq 1$, we first have
\[\begin{array}{l}
   n^{1-\frac{1}{\lambda}}\sum_{k\geq (nc)^{\frac{1}{\lambda}}} ck^{-\lambda}e^{-nck^{-\lambda}}
    =n^{1-\frac{1}{\lambda}}\sum_{k\geq 1} ck^{-\lambda}e^{-nck^{-\lambda}}1[k\geq (nc)^{\frac{1}{\lambda}}].
  \end{array}
\] Since $f_{n}(x)=n^{1-\frac{1}{\lambda}}ck^{-\lambda}e^{-nck^{-\lambda}}1[k\geq (nc)^{\frac{1}{\lambda}}]$  satisfies the condition of the Euler-Maclaurin Lemma with $x(n)=(nc)^{\frac{1}{\lambda}}$ and $f_{n}(x(n))\rightarrow 0$, we again have
\beq\begin{array}{l}
   n^{1-\frac{1}{\lambda}}\sum_{k\geq [(n+1)c]^{\frac{1}{\lambda}}} ck^{-\lambda}e^{-nck^{-\lambda}}
    = c\int_{1}^{\infty}n^{1-\frac{1}{\lambda}}x^{-\lambda}e^{-ncx^{-\lambda}}1[x\geq [(n+1)c]^{\frac{1}{\lambda}}]dx \\ \\
\hs = c\int_{[(n+1)c]^{\frac{1}{\lambda}}}^{\infty}n^{1-\frac{1}{\lambda}}x^{-\lambda}e^{-ncx^{-\lambda}}dx
    = c^{\frac{1}{\lambda}}\lambda^{-1}\Gamma\left(1-\frac{1}{\lambda}\right)\int_{0}^{(n+1)c} \frac{1}{\Gamma\left(1-\frac{1}{\lambda}\right)}s^{\left(1-\frac{1}{\lambda}\right)-1}e^{-s}ds \\ \\
\hs \rightarrow c^{\frac{1}{\lambda}}\lambda^{-1}\Gamma\left(1-\frac{1}{\lambda}\right)>0.
\end{array}\label{arg1}\eeq
On the other hand, for sufficiently large $n$, $I^{*}=\{k;p_{k}\leq \frac{1}{n+1}\} \subseteq \{k;k\geq k_{0}\}$, by parts 1) and 2) of Lemma \ref{2facts} at steps 2 and 4 below and (\ref{arg1}) at step 7, we have
 \[ \begin{array}{l}
    n^{1-1/\lambda}\sum_{k\in I^*}p_{k}(1-p_{k})^{n}\geq n^{1-1/\lambda}\sum_{k\in I^*}q_{k}(1-q_{k})^{n} \\ \\
    \hs \geq n^{1-1/\lambda}\sum_{k\in I^*}q_{k}\exp\left(-\frac{nq_{k}}{1-q_{k}}\right) \\ \\
    \hs \geq n^{1-1/\lambda}\sum_{k\in I^*}q_{k}\exp\left(-\frac{nq_{k}}{1-\sup_{I^*}q_{k}}\right)  \\ \\
    \hs \geq n^{1-1/\lambda}\sum_{k\in I^*}\exp(-2n(\sup_{I^*}q_{k})^{2})q_{k}e^{-nq_{k}} \\ \\
    \hs \geq n^{1-1/\lambda}\sum_{k\in I^*}\exp(-2/n)q_{k}e^{-nq_{k}} \\ \\
    \hs = \exp(-2/n)n^{1-1/\lambda}\sum_{k\in I^*}ck^{-\lambda}e^{-nck^{-\lambda}} \\ \\
    \hs \rightarrow c^{\frac{1}{\lambda}}\lambda^{-1}\Gamma\left(1-\frac{1}{\lambda}\right)>0.
 \end{array}
 \] Finally $t_{n}=n\sum_{k}p_{k}(1-p_{k})^{n}\geq n^{1/\lambda}n^{1-1/\lambda}\sum_{k\in I^*}p_{k}(1-p_{k})^{n}\rightarrow \infty$ as $n\rightarrow \infty$.
\hfill $\Box$

Theorem \ref{theorem2} puts distributions with power decaying tails, for example $p_{k}=c_{\lambda}k^{-\lambda}$, and those with slower decaying tails, for example $p_{k}=c_{\lambda} (k\ln^\lambda k)^{-1}$, where $\lambda>1$ and $c_{\lambda}>0$ is a constant which may depend on $\lambda$, in Domain 2.

In view of Lemma \ref{lemma1}, and Theorems \ref{theorem1} and \ref{theorem2}, Domain 1 has a more intuitive definition as given in the following lemma, the proof of which is trivial.
\begin{lemma} A distribution $P$ on $\mathscr{X}$ belongs to Domain 1 if and only if 1) the effective cardinality of $\mathscr{X}$ is $K=\infty$, and 2) $t_{n}\leq u_{P}$ for all $n$ and some constant $u_{P}>0$ which may depend on $P$.
\label{equivalence}
\end{lemma}

\begin{lemma}
 For any $P=\{p_{k}\}\in \mathscr{P}_{+}$, if there exists an integer $k_{0}\geq 1$ such that $p_{k}=c_{0}e^{-k}$ for all $k\geq k_{0}$ where $c_{0}>0$ is a constant, then
   \begin{enumerate}
     \item $t_{n}\leq u$  for some upper bound $u>0$; and
     \item $\lim_{n\rightarrow \infty}t_{n}$ does not exist.
   \end{enumerate}
\label{lemma5}
\end{lemma}

\ni {\it Proof}. Noting that the first finite terms of $t_{n}$ vanishes exponentially fast for any distribution, we may assume, without loss of generality, that $k_{0}=1$. For any given $n$, define $k^{*}=k^*(n)$ by
\beq \begin{array}{l}
    p_{k^*+1}< \frac{1}{n+1}\leq p_{k^*}.
    \label{defpk*1}
    \end{array}
\eeq
Noting that function $f_{n}(p)=np(1-p)^{n}$ increases for $p\in \left(0,\frac{1}{n+1}\right)$ decreases for $p\in \left(\frac{1}{n+1},1\right)$, we have for any $n$
\beq \begin{array}{rcll}
      f_{n}(p_{k})&\leq & f_{n}(p_{k^{*}}), & k\leq k^{*} \\ &&& \\
      f_{n}(p_{k})&<& f_{n}(p_{k^{*}}), & k\geq k^{*}+1.
     \end{array}
\label{fact1}
\eeq
Since $k^{*}=k^*(n)$ depends on $n$, we may express  $ p_{k^*}$ as, and define $c(n)$ by,
\beq \begin{array}{l}
   p_{k^*}=\frac{c(n)}{n}.
   \label{defpk*2}
   \end{array}
\eeq There are two main consequences of the expression in (\ref{defpk*2}). The first is that $t_{n}$ defined in (\ref{deftn}) may be expressed by (\ref{tn2}) below; and the second is that the sequence $c(n)$ perpetually oscillates between $1$ and $e$.

First, for each $n$, let us re-write each $p_{k}$ in terms of $p_{k^*}$, and therefore in terms of $n$ and $c(n)$.
\[ \begin{array}{c} p_{k^*+i}=e^{-i}\frac{c(n)}{n}\hspace{1em}\mbox{and} \hspace{1em}
   p_{k^*-j}=e^{j}\frac{c(n)}{n}\end{array}\]
for all appropriate positive integers $i$ and $j$. Therefore
\[  \begin{array}{rcl}
      f_{n}(p_{k^*+i})&=&ne^{-i}\frac{c(n)}{n}\left(1-e^{-i}\frac{c(n)}{n}\right)^{n}
        =\frac{c(n)}{e^{i}}\left(1-\frac{c(n)}{ne^{i}}\right)^{n}, \\ && \\
      f_{n}(p_{k^*-j})&=&ne^{j}\frac{c(n)}{n}\left(1-e^{j}\frac{c(n)}{n}\right)^{n}
         =c(n)e^{j}\left(1-\frac{c(n)e^{j}}{n}\right)^{n}.
   \end{array}
\] and
\vs
\beq \begin{array}{l}
   t_{n}=\sum_{k\leq k^*-1}f_{n}(p_{k})+f_{n}(p_{k^*})+\sum_{k\geq k^*+1}f_{n}(p_{k}) \\ \\
    \hs =c(n)\sum_{j=1}^{k^*-1}e^j\left(1-\frac{c(n)e^{j}}{n}\right)^{n}
      +c(n)\left(1-\frac{c(n)}{n}\right)^{n} +c(n)\sum_{i=1}^{\infty}e^{-i}\left(1-\frac{c(n)}{ne^{i}}\right)^{n}.
   \end{array}
   \label{tn2}
\eeq

\vs
Next we want to show that $c(n)$ oscillates perpetually over the interval $(n/(n+1),e)$ which approaches $[1,e)$ as $n$ increases indefinitely. This is so because, since $k^*$ is defined by (\ref{defpk*1}), we have
\[  \begin{array}{l}
    \frac{c(n)}{n}e^{-1}\leq \frac{1}{n+1}\leq \frac{c(n)}{n}
    \end{array}
\]or
\beq \begin{array}{l}
     e^{-1}<\frac{n}{n+1}\leq c(n)\leq \frac{n}{n+1}e<e.
     \label{specialINEQ}
     \end{array}
\eeq Furthermore by definition, $k^*=k^*(n)$  is an integer-valued increasing step function with unit increments. Let
$\{n_{k}; k\geq 1\}$ be the subsequence of $\mathbb{N}$ where $n_{k}$ is the positive integer value $n$ at which $k^*=k^*(n)$ jumps to a $k$ from $k-1$. Since
\[
  \begin{array}{c}
  c_{0}e^{-(k^*+1)}< \frac{1}{n+1}\leq c_{0}e^{-k^*} \\ \\
  e^{-(k^*+1)}< \frac{1}{c_{0}(n+1)}\leq e^{-k^*} \\ \\
  -(k^*+1)< -\ln (c_{0}(n+1))\leq -k^* \\ \\
  k^*+1 > \ln (c_{0}(n+1))\geq k^* ,
  \end{array}
\] we may write $k^*=\lfloor \ln(c_{0}(n+1)) \rfloor$ for each $n$. Clearly for each sufficiently large value $k^*$ there are multiple corresponding values of $n$ sharing the same value of $k^*$, denoted in the set $\{n_{k^*},n_{k^*}+1,\cdots,n_{k^*+1}-1\}$, and the size of the set increases indefinitely as $n\rightarrow \infty$.

Regarding the subsequence $\{n_{k^*}\}$ of $\mathbb{N}$, we have $1/n_{k^*}> p_{k^*}\geq 1/(n_{k^*}+1)$ or
\beq  \begin{array}{l}
       1-p_{k^*}\leq n_{k^*}p_{k^*}< 1,  \label{expression1}
       \end{array}
\eeq which implies that, for all sufficiently large $n$,

\beq \begin{array}{l}
      c(n_{k^*})=n_{k^*}p_{k^*}\in (1-\varepsilon,1) \label{expression2}
      \end{array}
\eeq where $\varepsilon>0$ is an arbitrarily small real value.


Similarly regarding the subsequence $\{n_{k^*+1}-1\}$ of $\mathbb{N}$, we first have
\[ \begin{array}{l}
     p_{k^*}=p_{k^*+1}e=\frac{n_{k^*+1}-1}{n_{k^*+1}-1}p_{k^*+1}e
       =\frac{1}{n_{k^*+1}-1}\left(\frac{n_{k^*+1}-1}{n_{k^*+1}}\right)(n_{k^*+1}p_{k^*+1})e
   \end{array}
\] and therefore by (\ref{expression1})
\[ \begin{array}{l}
     c(n_{k^*+1}-1)=\left(\frac{n_{k^*+1}-1}{n_{k^*+1}}\right)(n_{k^*+1}p_{k^*+1})e \rightarrow e
   \end{array}
\] which implies that, for all sufficiently large $n$,
\beq c(n_{k^*+1}-1)>e-\varepsilon
\label{expression3}
\eeq where $\varepsilon>0$ is an arbitrarily small real value. Furthermore over the set $\{n_{k^*},n_{k^*}+1,\cdots,n_{k^*+1}-1\}$, by the definition of $c(n)$ it is easy to see that $c(n)$ strictly increases with an exact increment of $p_{k^*}$ which decreases to zero as $n$ increases indefinitely. At this point, it has been established that the range of $c(n)$ for $n\geq n_{0}$, where $n_{0}$ is any positive integer, covers the entire interval $[1,e)$.

Noting $\mathbb{N}=\cup \{n_{k^*},n_{k^*}+1,\cdots,n_{k^*+1}-1\}$ where the union is over all possible integer values of $k^*$,
(\ref{expression2}) and (\ref{expression3}) jointly establish that the function $c(n)$ oscillates perpetually over the entire range of $[1,e)$.

The first part of the lemma follows that, noting that $e^{-1}\leq c(n)\leq e$ (see (\ref{specialINEQ})) and that $1-p\leq e^{-p}$ for all $p\in [0,1]$,
\[ \begin{array}{l}
   t_{n}=c(n)\sum_{j=1}^{k^*-1}e^j\left(1-\frac{c(n)e^{j}}{n}\right)^{n}
      +c(n)\left(1-\frac{c(n)}{n}\right)^{n} +c(n)\sum_{j=1}^{\infty}e^{-j}\left(1-\frac{c(n)}{ne^{j}}\right)^{n} \\ \\
   \hs \leq e\sum_{j=1}^{k^*-1}e^j\left(1-\frac{e^{j-1}}{n}\right)^{n}
      +e\left(1-\frac{e^{-1}}{n}\right)^{n} +e\sum_{j=1}^{\infty}e^{-j}\left(1-\frac{1}{ne^{j+1}}\right)^{n}  \\ \\
   \hs \leq e\sum_{j=1}^{k^*-1}e^je^{-e^{j-1}}
      +e\sum_{j=0}^{\infty}e^{-j}e^{-e^{-(j+1)}} \\ \\
   \hs \leq e^{2}\sum_{j=1}^{k^*-1}e^{j-1}e^{-e^{j-1}}
      +e^{2}\sum_{j=0}^{\infty}e^{-(j+1)}e^{-e^{-(j+1)}} \\ \\
   \hs < e^{2}\sum_{j=0}^{\infty}e^je^{-e^{j}}
      +e^{2}\sum_{j=1}^{\infty}e^{-j}e^{-e^{-j}}\vcentcolon=u.
   \end{array}
\]

\vs
For the second part of the lemma, consider, for any fixed $c>0$,
\[ \begin{array}{l}
     t^{*}_{n}=c\sum_{j=1}^{k^*-1}e^j\left(1-\frac{ce^{j}}{n}\right)^{n}
      +c\left(1-\frac{c}{n}\right)^{n} +c\sum_{j=1}^{\infty}e^{-j}\left(1-\frac{c}{ne^{j}}\right)^{n}.
   \end{array}
\]By Dominated Convergence Theorem,
\[ \begin{array}{l}
   t(c)\vcentcolon=\lim_{n\rightarrow \infty} t^{*}_{n}=c\sum_{j=0}^{\infty}e^je^{-ce^{j}}+c\sum_{j=1}^{\infty}e^{-j}e^{-ce^{-j}},
   \end{array}
\] and $t(c)$ is a non-constant function in $c$ on $[1,e]$.

The argument thus far implies that, as $n$ increases, $c(n)$ repeatedly visits any arbitrarily small closed interval $[a,b]\subset [1,e]$ infinitely often, and therefore there exists for each such interval a subsequence $\{n_{l}; l\geq 1\}$ of $\mathbb{N}$ such that $c(n_{l})$ converges, {\it i.e.}, $c(n_{l})\rightarrow \theta$ for some $\theta \in [a,b]$. Since $t(c)$ is a non-constant function on $[1,e]$, there exist two non-overlapping closed intervals, $[a_{1},b_{1}]$ and $[a_{2},b_{2}]$ in $[1,e]$, satisfying
 \[ \begin{array}{l}
    \max_{a_{1}\leq c\leq b_{1}}t(c)<\min_{a_{2}\leq c\leq b_{2}}t(c),
    \end{array}
 \] such that there exist two sub-sequences of $\mathbb{N}$, said $\{n_{l}; l\geq 1\}$ and $\{n_{m}; m\geq 1\}$, such that $c(n_{l})\rightarrow \theta_{1}$ for some $\theta_{1} \in [a_{1},b_{1}]$ and $c(n_{m})\rightarrow \theta_{2}$ for some $\theta_{2} \in [a_{2},b_{2}]$.

Consider the limit of $t_{n}$ along $\{n_{l}; l\geq 1\}$, again by Dominated Convergence Theorem,
\[\begin{array}{l}
  \lim_{n_{l}\rightarrow \infty} t_{n_{l}}=\lim_{n_{l}\rightarrow \infty}\left[c(n_{l})\sum_{j=0}^{k^*-1}e^j\left(1-\frac{c(n_{l})e^{j}}{n}\right)^{n}
      +c(n_{l})\sum_{j=1}^{\infty}e^{-j}\left(1-\frac{c(n_{l})}{ne^{j}}\right)^{n}\right] \\ \\
  \hs = \theta_{1}\sum_{j=0}^{\infty}e^je^{-\theta_{1}e^{j}}+\theta_{1}\sum_{j=1}^{\infty}e^{-j}e^{-\theta_{1}e^{-j}}
      = t(\theta_{1}).
\end{array}\] A similar argument gives $\lim_{n_{m}\rightarrow \infty} t_{n_{m}}=t(\theta_{2})$, but $t(\theta_{1})\neq t(\theta_{2})$ by construction, and hence $\lim_{n\rightarrow\infty}t_{n}$ does not exist.
\hfill $\Box$

A similar proof to that of Lemma \ref{lemma5} immediately gives Theorem \ref{theorem3} below with a slightly more general statement.
\begin{theorem}
 For any given probability distribution $P=\{p_{k}; k\geq 1\}$, if there exists constants $a>1$ and integer $k_{0}\geq 1$ such that for all $k\geq k_{0}$
  \beq p_{k}= ca^{-k},\label{conditionT3}\eeq then
   \begin{enumerate}
     \item $t_{n}\leq u_{a}$  for some upper bound $u_{a}>0$ which may depend on $a$; and
     \item $\lim_{n\rightarrow \infty}t_{n}$ does not exist.
   \end{enumerate}
\label{theorem3}
\end{theorem}
\vs
Theorem \ref{theorem3} puts distributions with tails of geometric progression, for example $p_{k}=c_{\lambda}e^{-\lambda k}$ where $\lambda >0$ and $c_{\lambda}>0$ are constants or $p_{k}=2^{-k}$, in Domain 1.

Next we develop a notion of relative dominance of one probability distribution over another on a countable alphabet within $\mathscr{P}_{+}$. Let $\#A$ denote the cardinality of a set $A$.
\begin{definition} Let $Q^*\in \mathscr{P}_{+}$ and $P\in \mathscr{P}_{+}$ be two distributions on $\mathscr{X}$, and let $Q=\{q_{k}\}$ be a non-increasingly ordered version of $Q^*$.
 $Q^*$ is said to dominate $P$ if
 \[
    \#\{i; p_{i}\in (q_{k+1},q_{k}], i\geq 1 \}\leq M <\infty
 \] for every $k\geq 1$, where $M$ is a finite positive integer.
\label{defRegularity}
\end{definition}

It is easy to see that the notion of dominance by Definition \ref{defRegularity} is a tail property, and that
it is transitive, {\it i.e.}, if $P_{1}$ dominates $P_{2}$ and $P_{2}$ dominates $P_{3}$ then $P_{1}$ dominates $P_{3}$. It says in essence that, if $P$ is dominated by $Q$, then the $p_{i}$s do not get overly congregated locally into some intervals defined by the $q_{k}$s.

The following examples shed a bit of intuitive light on the notion of dominance by Definition \ref{defRegularity}.

\begin{example}Let $p_{k}=c_{1}e^{-k^{2}}$ and $q_{k}=c_{2}e^{-k}$ for all $k\geq k_{0}$ for some integer $k_{0}\geq 1$ and other two constants $c_{1}>0$ and $c_{2}>0$. For every sufficiently large $k$, suppose
$p_{j}=c_{1}e^{-j^{2}}\leq q_{k}=c_{2}e^{-k}$, then $-j^{2}\leq \ln\left(c_{2}/c_{1}\right)-k$ and
$j+1\geq [k+\ln\left(c_{1}/c_{2}\right)]^{1/2}+1$. It follows that
\[
\begin{array}{l}
     p_{j+1}=c_{1}e^{-(j+1)^2}\leq c_{1}e^{-\left(\sqrt{k+\ln\left(c_{1}/c_{2}\right)}+1\right)^2}
         =c_{1}e^{-\left(k+\ln\left(c_{1}/c_{2}\right)+1\right)-2\sqrt{k+\ln\left(c_{1}/c_{2}\right)}} \\ \\
       \hs  =c_{2}e^{-\left(k+1\right)-2\sqrt{k+\ln\left(c_{1}/c_{2}\right)}}
         =c_{2}e^{-\left(k+1\right)}e^{-2\sqrt{k+\ln\left(c_{1}/c_{2}\right)}}\leq c_{2}e^{-\left(k+1\right)}=q_{k+1}.
\end{array}
\] This means that if $p_{j}\in (q_{k+1},q_{k}]$ then necessarily $p_{j+1}\not \in (q_{k+1},q_{k}]$, which implies that each
interval $(q_{k+1},q_{k}]$ can contain only one $p_{j}$ at most for a sufficiently large $k$, {\it i.e.}, $k\geq k_{00}\vcentcolon=\max\{k_{0}, \ln(c_{2}/c_{1})\}$. Since there are only finite $p_{j}$s covered by $\cup_{1\leq k<k_{00}}(q_{k},q_{k+1}]$,
$Q=\{q_{k}\}$ dominates $P=\{p_{i}\}$.
\label{example1}
\end{example}

\begin{example}Let $p_{k}=c_{1}a^{-k}$ and $q_{k}=c_{2}b^{-k}$ for all $k\geq k_{0}$ for some integer $k_{0}\geq 1$ and other two constants $a>b>1$. For every sufficiently large $k$, suppose
$p_{j}=c_{1}a^{-j}\leq q_{k}=c_{2}b^{-k}$, then
$-j\ln a\leq \ln\left(c_{2}/c_{1}\right)-k\ln b$ and $j+1\geq k(\ln b/\ln a)+1+\ln\left(c_{1}/c_{2}\right)/\ln a$. It follows that
\[ \begin{array}{l}
     p_{j+1}=c_{1}a^{-\left( k\frac{\ln b}{\ln a}+1+\frac{\ln\left(c_{1}/c_{2}\right)}{\ln a} \right)  }
            =c_{1}a^{-\left( k\log_{a}{b}+1+\frac{\ln\left(c_{1}/c_{2}\right)}{\ln a} \right)  } \\ \\
     \hs =c_{1}b^{-k}a^{-1} a^{ -\frac{\ln\left(c_{1}/c_{2}\right)}{\ln a}}
         \leq c_{1}b^{-(k+1)} a^{ -\log_{a}\left(c_{1}/c_{2}\right)}=c_{2}b^{-(k+1)}=q_{k+1}.
\end{array}
\] By a similar argument as that in Example \ref{example1},
$Q=\{q_{k}\}$ dominates $P=\{p_{i}\}$.
\label{example2}
\end{example}

\begin{example}
 Let $p_{k}=c_{1}k^{-r}e^{-\lambda k}$ for some integer $k_{0}\geq 1$ and constants $\lambda>0$ and $r>0$, and $q_{k}=c_{2}e^{-\lambda k}$ for all $k\geq k_{0}$. Suppose for a $k\geq k_{0}$ there is a $j$ such that $p_{j}=c_{1}j^{-r}e^{-\lambda j}\in (q_{k+1}=c_{2}e^{-\lambda (k+1)},q_{k}=c_{2}e^{-\lambda k}]$, then
 \[ \begin{array}{l}
    p_{j+1}=c_{1}(j+1)^{-r}e^{-\lambda (j+1)}=c_{1}(j+1)^{-r}e^{-\lambda j}e^{-\lambda}
    \leq c_{1}j^{-r}e^{-\lambda j}e^{-\lambda} \\ \\
    \hs \leq c_{2}e^{-\lambda k}e^{-\lambda}=q_{k+1},
    \end{array}
\] which implies that there is at most one $p_{j}$ in $(q_{k+1},q_{k}]$ for every sufficiently large $k$. Therefore $Q=\{q_{k}\}$ dominates $P=\{p_{i}\}$.
\label{example3}
\end{example}

\begin{example}
 Let $p_{k}=c_{1}k^{r}e^{-\lambda k}$ for some integer $k_{0}\geq 1$ and constants $\lambda>0$ and $r>0$, and $q_{k}=c_{2}e^{-(\lambda/2) k}$ for all $k\geq k_{0}$. Suppose for any sufficiently large $j$, $j\geq j_{0}\vcentcolon=\left[e^{\lambda/(2r)}-1\right]^{-1}$, we have $p_{j}=c_{1}j^{r}e^{-\lambda j}\in \left(q_{k+1}=c_{2}e^{-(\lambda/2) (k+1)},q_{k}=c_{2}e^{-(\lambda/2) k}\right]$ for some sufficiently large $k\geq k_{0}$, then
 \[ \begin{array}{l}
    p_{j+1}=c_{1}(j+1)^{r}e^{-\lambda (j+1)}=c_{1}(j+1)^{r}e^{-\lambda j}e^{-\lambda}
       =c_{1}j^{r}e^{-\lambda j}e^{-\lambda}\frac{(j+1)^{r}}{j^{r}}  \\ \\
    \hs \leq c_{2}e^{-\frac{\lambda}{2} k}e^{-\lambda} \left(\frac{j+1}{j}\right)^{r}
    = c_{2}e^{-\frac{\lambda}{2} (k+1)}e^{-\frac{\lambda}{2}} \left(\frac{j+1}{j}\right)^{r}   \\ \\
    \hs \leq q_{k+1}e^{-\frac{\lambda}{2}} \left(\frac{j_{0}+1}{j_{0}}\right)^{r}=q_{k+1}
    \end{array}
\] which implies that there is at most one $p_{j}$ in $(q_{k+1},q_{k}]$ for every sufficiently large $k$. Therefore $Q=\{q_{k}\}$ dominates $P=\{p_{i}\}$.
\label{example4}
\end{example}

\begin{example}Let $p_{k}=q_{k}$ for all $k\geq 1$. $Q=\{q_{k}\}$ and $P=\{p_{k}\}$ dominate each other.
\label{example5}
\end{example}

While in each of Examples \ref{example1} through \ref{example4}, the dominating distribution $Q$ has a thicker tail than $P$ in the usual sense, the dominance of Definition \ref{defRegularity} in general is not implied by such a thinner/thicker tail relationship. This is so because a distribution $P\in \mathscr{P}_{+}$, satisfying $p_{k}\leq q_{k}$ for all sufficiently large $k$, could exist yet congregate irregularly to have an unbounded $\sup_{k\geq 1}\#\{p_{i};p_{i}\in (q_{k+1},q_{k}], i\geq 1\}$. One such example is given in Section 3 below. In this regard, the dominance of Definition \ref{defRegularity} is more appropriately considered as a regularity condition. However it may be interesting to note that the said regularity is a relative one in the sense that the behavior of $P$ is regulated by a reference distribution $Q$. This relative regularity gives an
umbrella structure in Domain 1 as demonstrated by the theorem below.

\begin{theorem} If two distributions $P$ and $Q$ in $\mathscr{P}_{+}$ on a same countably infinite alphabet $\mathscr{X}$ are such that $Q$ is in Domain 1 and $P$ is dominated by $Q$, then $P$ belongs to Domain 1.
\label{theorem4}
\end{theorem}
\ni {\it Proof}. Without loss of generality, it may be assumed that $Q$ is non-increasingly ordered. For every $n$, there exists a $k_{n}$ such that $\frac{1}{n+1}\in (q_{k_{n}+1},q_{k_{n}}]$. Noting that
the function $np(1-p)^{n}$ increases in $p$ over $(0,1/(n+1)]$, attains its maximum value of $[1-1/(n+1)]^{n+1}<e^{-1}$ at $p=1/(n+1)$, and decreases over $[1/(n+1), 1]$, consider
\[\begin{array}{l}
  t_{n}(P)=\sum_{k\geq 1}np_{k}(1-p_{k})^{n} \\ \\
  \hs  = \sum_{k; p_{k}\leq q_{k_{n}+1}}np_{k}(1-p_{k})^{n}+\sum_{k; q_{k_{n}+1}< p_{k}\leq q_{k_{n}}}np_{k}(1-p_{k})^{n}
        +\sum_{k; p_{k}> q_{k_{n}}}np_{k}(1-p_{k})^{n} \\ \\
  \hs \leq M\sum_{k\geq k_{n}+1}nq_{k}(1-q_{k})^{n}+\sum_{k; q_{k_{n}+1}< p_{k}\leq q_{k_{n}}}e^{-1}
        +M\sum_{1\leq k\leq k_{n}}nq_{k}(1-q_{k})^{n} \\ \\
  \hs = M\sum_{k\geq 1}nq_{k}(1-q_{k})^{n}+\sum_{k; q_{k_{n}+1}< p_{k}\leq q_{k_{n}}}e^{-1} \\ \\
  \hs \leq M t_{n}(Q)+Me^{-1} < \infty.
 \end{array}
\] The desired result immediately follows.
\hfill $\Box$

\begin{corollary} Any distribution $P$ on a countably infinite alphabet $\mathscr{X}$ satisfying $p_{k}=ae^{-\lambda k}$, $p_{k}=be^{-\lambda k^2}$,
or $p_{k}=ck^{r}e^{-\lambda k}$ for all $k\geq k_{0}$, where $k_{0}\geq 1$, $\lambda>0$, $r\in (-\infty,+\infty)$, $a>0$, $b>0$ and $c>0$ are constants, is in Domain 1.
\label{coro2}
\end{corollary}
\ni {\it Proof}. The result is immediate following Theorem \ref{theorem4} and Examples \ref{example1} through \ref{example4}.   \hfill $\Box$

\section{Constructed Examples.}
The first constructed example shows that the notion of thinner tail, in the sense of $p_{k}\leq q_{k}$ for $k\geq k_{0}$ where $k_{0}\geq 1$ is some fixed integer and $P=\{p_{k}\}$ and $Q=\{q_{k}\}$ are two distributions, does not imply a dominance of $Q$ over $P$.
\begin{example} Consider any strictly decreasing distribution $Q=\{q_{k};k\geq 1\}\in \mathscr{P}_{+}$ and the following grouping of the index set $\{k;k\geq 1\}$.
\[ \begin{array}{l}
G_{1}=\{1\},G_{2}=\{2,3\},\cdots, G_{m}=\{m(m-1)/2+1,\cdots,m(m-1)/2+m\},\cdots.
\end{array}
\] $\{G_{m};m\geq 1\}$ is a partition of the index set $\{k;k\geq 1\}$ and each group $G_{m}$ contains $m$ consecutive indices. A new distribution $P=\{p_{k}\}$ is constructed according to the following steps:
\begin{enumerate}
  \item For each $m\geq 2$, let $p_{k}=q_{m(m-1)/2+m}$ for all $k\in G_{m}$.
  \item $p_{1}=1-\sum_{k\geq 2}p_{k}$.
\end{enumerate}
In the first step, $m(m-1)/2+m=m(m+1)/2$ is the largest index in $G_{m}$ and therefore $q_{m(m+1)/2}$ is the smallest $q_{k}$ with index $k\in G_{m}$. Since
\[ \begin{array}{l}
    0\leq \sum_{k\geq 2}p_{k}=\sum_{m\geq 2}mq_{m(m+1)/2}< \sum_{k\geq 2}q_{k}\leq 1,
\end{array}\] $p_{1}$ so assigned is a probability. The distribution $P=\{p_{k}\}$ satisfies $p_{k}\leq q_{k}$ for every $k\geq 2=k_{0}$. However the number of terms of $p_{i}$ in the interval $(q_{m(m+1)/2+1},q_{m(m+1)/2}]$ is at least $m$ and it increases indefinitely as $m\rightarrow \infty$; and hence $Q$ does not dominate $P$.
\label{ex6}
\end{example}

The second constructed example shows that the notion of the dominance of $Q=\{q_{k}\}$ over $P=\{p_{k}\}$, as defined in Definition \ref{defRegularity}, does not imply that $P$ has thinner tail than $Q$, in the sense of $p_{k}\leq q_{k}$ for $k\geq k_{0}$ where $k_{0}\geq 1$ is some fixed integer.
\begin{example} Consider any strictly decreasing distribution $Q=\{q_{k};k\geq 1\}\in \mathscr{P}_{+}$ and the following grouping of the index set $\{k;k\geq 1\}$.
\[ \begin{array}{l}
G_{1}=\{1,2\},G_{2}=\{3,4\},\cdots, G_{m}=\{2m-1,2m\},\cdots.
\end{array}
\] $\{G_{m};m\geq 1\}$ is a partition of the index set $\{k;k\geq 1\}$ and each group $G_{m}$ contains $2$ consecutive indices, the first one odd and the second one even. The construction of a new distribution $P=\{p_{k}\}$ is as follows: for each group $G_{m}$ with its two indices $k=2m-1$ and $k+1=2m$,
let $p_{k}=p_{k+1}=(q_{k}+q_{k+1})/2$. With the new distribution $P=\{p_{k}\}$ so defined, we have $p_{2m}< q_{2m}$ and $p_{2m-1}>q_{2m-1}$ for all $m\geq 1$. Clearly $Q$ dominates $P$ ($P$ dominates $Q$ as well), but $P$ does not have a thinner tail in the usual sense.
\end{example}

At this point, it becomes clear that the notation of dominance of Definition \ref{defRegularity} and the notation of thinner/thicker tail in the usual sense are two independent notions.

The next constructed example below shows that there exists a distribution such that he associated $t_{n}$ approaches infinity along one subsequence of $n$ and is bounded above along another subsequence of $n$, hence belonging to Domain $T$. Domain $T$ is not empty.

\begin{example}  Consider the probability sequence $q_{j}=2^{-j}$, for $j=1,2,\cdots$, along with a diffusion sequence $d_{i}=2^{i}$, for $i=1,2,\cdots$. A probability sequence $\{p_{k}\}$, for $k=1,2,\cdots$, is constructed by the following steps:
\begin{enumerate}
   \item[\hs $1^{st}\vcentcolon$] \begin{enumerate}
       \item Take the first value of $d_{i}$, $d_{1}=2^1$, and assign the first $2d_{1}=2^2=4$ terms of $q_{j}$,
           $q_{1}=2^{-1}, q_{2}=2^{-2}, q_{3}=2^{-3}, q_{4}=2^{-4}$, to the first $4$ terms of $p_{k}$,
           $p_{1}=2^{-1}, p_{2}=2^{-2}, p_{3}=2^{-3},p_{4}=2^{-4}$.
       \item Take the next unassigned term in $q_{j}$, $q_{5}=2^{-5}$, and diffuse it into $d_{1}=2$ equal terms, $2^{-6}$ and $2^{-6}$.
          \begin{enumerate}
          \item Starting at $q_{5}$ in the sequence $\{q_{j}\}$, look forwardly ($j>5$) for terms greater or equal to $2^{-6}$, if any, continue to assign them to $p_{k}$. In this case, there is only one such term $q_{6}=2^{-6}$ and it is assigned to $p_{5}=2^{-6}$.
          \item Take the $d_{1}=2$ diffused terms and assign them to $p_{6}=2^{-6}$ and $p_{7}=2^{-6}$. At this point, the first few terms of the partially assigned sequence $\{p_{k}\}$ are
            \[ p_{1}=2^{-1}, p_{2}=2^{-2}, p_{3}=2^{-3},p_{4}=2^{-4},p_{5}=2^{-6},p_{6}=2^{-6}, p_{7}=2^{-6}.\]
          \end{enumerate}
       \end{enumerate}
  \item[\hs $2^{nd}\vcentcolon$] \begin{enumerate}
       \item Take the next value of $d_{i}$, $d_{2}=2^2$, and assign the next $2d_{2}=2^3=8$ unused terms of $q_{j}$,
           $q_{7}=2^{-7}, \cdots, q_{14}=2^{-14}$, to the next $8$ terms of $p_{k}$,
           $p_{8}=2^{-7}, \cdots, p_{15}=2^{-14}$.
       \item Take the next unassigned term in $q_{j}$, $q_{15}=2^{-15}$, and diffuse it into $d_{2}=4$ equal terms of $2^{-17}$ each.
          \begin{enumerate}
            \item Starting at $q_{15}$ in the sequence of $\{q_{j}\}$, look forwardly ($j>15$) for terms greater or equal to $2^{-17}$, if any, continue to assign them to $p_{k}$. In this case, there are 2 such terms $q_{16}=2^{-16}$ and $q_{17}=2^{-17}$, and they are assigned to $p_{16}=2^{-16}$ and $p_{17}=2^{-17}$.
            \item Take the $d_{2}=2^2=4$ diffused terms and assign them to $p_{18}=2^{-17}, \cdots, p_{21}=2^{-17}$. At this point, the first few terms of the partially assigned sequence $\{p_{k}\}$ are
            \[ \begin{array}{c}p_{1}=2^{-1}, p_{2}=2^{-2}, p_{3}=2^{-3},p_{4}=2^{-4},\\ \\
            p_{5}=2^{-6},p_{6}=2^{-6}, p_{7}=2^{-6}, \\ \\
            p_{8}=2^{-7}, p_{9}=2^{-8},\cdots, p_{15}=2^{-14}, p_{16}=2^{-16}, \\ \\
            p_{17}=2^{-17}, p_{18}=2^{-17}, \cdots, p_{21}=2^{-17}.
            \end{array}\]
          \end{enumerate}
       \end{enumerate}
  \item[\hs $i^{th}\vcentcolon$] \begin{enumerate}
       \item In general, take the next value of $d_{i}$, say $d_{i}=2^{i}$, and assign the next $2d_{i}=2^{i+1}$ unused terms of $q_{j}$, say
           $q_{j_{0}}=2^{-j_{0}}, \cdots, q_{j_{0}+2^{i+1}-1}=2^{-(j_{0}+2^{i_{0}+1}-1)}$, to the next $2d_{i}=2^{i+1}$ terms of $p_{k}$, say $p_{k_{0}}=2^{-j_{0}}, \cdots, p_{k_{0}+2^{i+1}-1}=2^{-(j_{0}+2^{i+1}-1)}$.
       \item Take the next unassigned term in $q_{j}$, $q_{j_{0}+2^{i+1}}=2^{-(j_{0}+2^{i+1})}$, and diffuse it into $d_{i}=2^{i}$ equal terms, $2^{-(j_{0}+i+2^{i+1})}$ each.
            \begin{enumerate}
              \item Starting at $q_{j_{0}+2^{i+1}}$ in the sequence of $\{q_{j}\}$, look forwardly ($j>j_{0}+2^{i+1}$) for terms greater or equal to $2^{-(j_{0}+i+2^{i+1})}$, if any, continue to assign them to $p_{k}$. Denote the last assigned $p_{k}$ as $p_{k_{0}}$.
              \item Take the $d_{i}=2^{i}$ diffused terms and assign them to $p_{k_{0}+1}=2^{-(j_{0}+i+2^{i+1})}$, $\cdots$, $p_{k_{0}+2^{i}}=2^{-(j_{0}+i+2^{i+1})}$.
            \end{enumerate}
       \end{enumerate}
\end{enumerate}

In essence, the sequence $\{p_{k}\}$ is generated based on the sequence $\{q_{j}\}$ with infinitely many selected $j$'s at each of which $q_{j}$ is diffused into increasingly many equal probability terms according a diffusion sequence $\{d_{i}\}$. The diffused sequence is then re-arranged in a non-increasing order. By construction, it is clear that the sequence $\{p_{k}; k\geq 1\}$, satisfies the following properties:
\begin{enumerate}
 \item[${\cal A}_{1}$:] $\{p_{k}\}$ is a probability sequence in a non-increasing order.
 \item[${\cal A}_{2}$:] As $k$ increases, $\{p_{k}\}$ is a string of segments alternating between two different types: 1) a strictly decreasing segment and 2) a segment (a run) of equal probabilities.
 \item[${\cal A}_{3}$:] As $k$ increases, the length of the last run increases and approaches infinity.
 \item[${\cal A}_{4}$:] In each run, there are exactly $d_{i}+1$ equal terms, $d_{i}$ of which are diffused terms and 1 of which belongs to the original sequence $q_{j}$.
 \item[${\cal A}_{5}$:] Between two consecutive runs (with lengths $d_{i}+1$ and $d_{i+1}+1$ respectively), the strictly decreasing segment in the middle has at least $2 d_{i+1}=4d_{i}= d_{i}+3d_{i}>d_{i}+d_{i+1}$ terms.
 \item[${\cal A}_{6}$:] For any $k$, $1/p_{k}$ is a positive integer.
\end{enumerate}

Next we want to show that there is a subsequence $\{n_{i}\} \in \mathbb{N}$ such that $t_{n_{i}}$ defined with $\{p_{k}\}$ approaches infinity. Toward that end, consider the subsequence $\{p_{k_{i}}; i\geq 1\}$ of $\{p_{k}\}$ where the index $k_{i}$ is such that $p_{k_{i}}$ is first term in the $i^{th}$ run segment. Let $\{n_{i}\}=\{1/p_{k_{i}}\}$ which by
${\cal A}_{6}$ is a subsequence of $\mathbb{N}$. By ${\cal A}_{3}$ and ${\cal A}_{4}$,
\[ \begin{array}{l}
  t_{n_{i}}=n_{i}\sum_{k\geq 1}p_{k}(1-p_{k})^{n_{i}}>n_{i}(d_{i}+1)p_{k_{i}}(1-p_{k_{i}})^{n_{i}}
  =(d_{i}+1)\left(1-\frac{1}{n_{i}}\right)^{n_{i}}\rightarrow \infty.
\end{array}
\]
Consider next the subsequence $\{p_{k_{i}-(d_{i}+1)};i\geq 1\}$ of $\{p_{k}\}$ where the index $k_{i}$ is such that $p_{k_{i}}$ is first term in the $i^{th}$ run segment, and therefore $p_{k_{i}-(d_{i}+1)}$ is the $(d_{i}+1)^{th}$ term counting backwards from $p_{k_{i}-1}$,
into the preceding segment of at least $2d_{i}$ strictly decreasing terms. Let $\{m_{i}\}=\{1/p_{k_{i}-(d_{i}+1)}-1\}$ (so
$p_{k_{i}-(d_{i}+1)}=(m_{i}+1)^{-1}$) which by
${\cal A}_{6}$ is a subsequence of $\mathbb{N}$.
\[ \begin{array}{l}
  t_{m_{i}}=m_{i}\sum_{k\geq 1}p_{k}(1-p_{k})^{m_{i}}=m_{i}\sum_{k\leq k_{i}-(d_{i}+1)}p_{k}(1-p_{k})^{m_{i}}+m_{i}\sum_{k\geq k_{i}-d_{i}}p_{k}(1-p_{k})^{m_{i}} \\ \\
  \hs \vcentcolon=t_{m_{i},1}+t_{m_{i},2}.
\end{array}
\]  Before proceeding further, let us note several detailed facts. First, the function $np(1-p)^{n}$ increases in $[0,1(n+1)]$, attains maximum at $p=1/(n+1)$, and decreases in $[1/(n+1),1]$. Second, since $p_{k_{i}-(d_{i}+1)}=(m_{i}+1)^{-1}$,
by ${\cal A}_{1}$ each summand in $t_{m_{i},1}$ is bounded above by  $m_{i}p_{k_{i}-(d_{i}+1)}(1-p_{k_{i}-(d_{i}+1)})^{m_{i}}$ and each
summand in $t_{m_{i},2}$ is bounded above by $m_{i}p_{k_{i}-d_{i}}(1-p_{k_{i}-d_{i}})^{m_{i}}$. Third, by ${\cal A}_{4}$ and ${\cal A}_{5}$, for each diffused term of $p_{k'}$ with $k'\leq k_{i}-(d_{i}+1)$ in a run there is a different non-diffused term $p_{k''}$ with $k''\leq k_{i}-(d_{i}+1)$ such that  $p_{k'}>p_{k''}$ and therefore $m_{i}p_{k'}(1-p_{k'})^{m_{i}}\leq m_{i}p_{k''}(1-p_{k''})^{m_{i}}$; and similarly, for each diffused term of $p_{k'}$ with $k'\geq k_{i}-d_{i}$ in a run there is a different non-diffused term $p_{k''}$ with $k''\geq k_{i}-d_{i}$ such that  $p_{k'}<p_{k''}$ and therefore $m_{i}p_{k'}(1-p_{k'})^{m_{i}}\leq m_{i}p_{k''}(1-p_{k''})^{m_{i}}$. These facts imply that
\[\begin{array}{l}
  t_{m_{i}}=t_{m_{i},1}+t_{m_{i},2}=m_{i}\sum_{k\leq k_{i}-(d_{i}+1)}p_{k}(1-p_{k})^{m_{i}}+m_{i}\sum_{k\geq k_{i}-d_{i}}p_{k}(1-p_{k})^{m_{i}} \\ \\
  \hs \leq 2m_{i}\sum_{j\geq 1}q_{j}(1-q_{j})^{m_{i}}< \infty
\end{array}
\] and the last inequality above is due to Corollary \ref{coro2}.
\end{example}

\section{A Statistical Implication.}
While the domains of attraction on alphabets have probabilistic merit, the statistical implication is also quite significant. Zhang and Zhou (2010) showed that $\zeta_{1,v}$ is estimable (there exists at least one unbiased estimator of $\zeta_{1,v}$), and established an unbiased estimator of $\zeta_{1,v}$ for every $v\leq n-1$. Their estimator is
 \beq \begin{array}{l}\label{z1v}
 Z_{1,v}=\frac{n^{1+v}[n-(1+v)]!}{n!}\sum_{k\geq 1}
   \left[ \hat{p}_{k} \prod_{j=0}^{v-1}\left(1-\hat{p}_{k}-\frac{j}{n}\right)\right].
\end{array}\eeq
Therefore there readily exists an unbiased estimator of $t_{v}$ for every $v\leq n-1$ namely
 \beq \begin{array}{l}\label{thatv}
 \hat{t}_{v}=vZ_{1,v}.
\end{array}\eeq Zhang and Zhou (2010) also established several useful statistical properties of $\hat{t}_{v}$, including the asymptotic normality and that $\hat{t}_{v}$ is the uniformly minimum variance unbiased estimator ($umvue$) when $K<\infty$.

The availability of $\hat{t}_{v}$ gives much added merit to the discussion of the domains of attraction on alphabets as presented in this paper. Specifically the fact that the asymptotic behavior of $t_{n}$ characterizes the tail probability of the underlying $P$ and the fact that the trajectory of $t_{v}$ up to $v=n-1$ is estimable suggest that much could be revealed by a sufficiently large sample.

\end{document}